\documentclass[leqno]{amsart} 

\newcommand{\R}{{\mathbb R}}

\newcommand {\gd}{\displaystyle}

\newcommand {\rt}{\rightarrow}

\newcommand {\fraco}{\rightharpoonup}

\begin{document}

\title[Euler equations in noncylindrical domains]
{Two-dimensional incompressible ideal flows in a noncylindrical material domain}

\author{F. Z. Fernandes, M. C. Lopes Filho}
 
%\address{Fl\'avia Zechineli Fernandes\hfill\break
%Departamento de Matemática,
%Universidade Estadual de Campinas,
%Campinas, SP 13082-600, Brasil, and \hfill\break 
%Milton da Costa Lopes Filho\hfill\break 
%Departamento de Matemática,
%Universidade Estadual de Campinas, Campinas, SP 13082-600, Brasil} 
%\email{zechinel@gmail.com, mlopes@ime.unicamp.br}

%\date{}
%$\thanks{Submitted December 16, 2003. Published February 25, 2004.}
%\subjclass[2000]{34B05, 34B24, 92D25} 
%\keywords{Anti,
%\hfill\break\indent
%nonlinear boundary value problem}

\begin{abstract} The purpose of this work is to prove existence of a weak solution of the 
two dimensional incompressible Euler equations on a noncylindrical domain consisting of a 
smooth, bounded, connected and simply connected domain undergoing a prescribed motion. 
We prove existence of a weak solution for initial vorticity in $L^p$, for $p>1$. This
work complements a similar result by C. He and L. Hsiao, who proved existence assuming 
that the flow velocity is tangent to the moving boundary, see [JDE v. 163 (2000) 265--291].
\end{abstract}

\maketitle
\numberwithin{equation}{section}
\newtheorem{theorem}{Theorem}[section]
\newtheorem{lemma}[theorem]{Lemma}
\newtheorem{remark}[theorem]{Remark}
\newtheorem{proposition}[theorem]{Proposition}
\newtheorem{definition}[theorem]{Definition}

%--------------------------------------------------------------------------------------------------------------------
\section{Introduction}

In this work we will prove the existence of weak solutions of the incompressible Euler equations in two-dimensional 
domains with smoothly moving boundaries.  Previous work on incompressible flow in noncylindrical domains has 
addressed both viscous and ideal flow, and both weak and strong solutions. Existence of weak solutions for the 
Navier-Stokes equations on a noncylindrical domain was first studied by H. Fujita and N. Sauer, see cite{fujita}, 
whose work was later complemented by D. Bock in \cite{bock}, A. Inoue and M. Wakimoto in \cite{inoue} and T. 
Miyakawa and Y. Teramoto in \cite{teramoto}. For ideal flow, well-posedness of the two-dimensional problem in 
the case of smooth solutions was studied by H. Kozono in \cite{kozono} and the existence of a weak solution was 
studied by C. He and L. Hsiao in \cite{he}.  In their work, He and Hsiao assumed that the flow velocity is tangent 
to the moving boundary for each fixed time. As a consequence, their result does not include existence of a weak 
solution for the case of a noncylindrical material domain with prescribed motion. 
The purpose of this paper is to complement He and Hsiao's argument to prove existence of weak solutions for the incompressible two-dimensional Euler equations in a noncylindrical domain which moves with the flow 
in a smoothly prescribed manner.

We describe our treatment of this problem as follows. We start from a prescribed movement of the domain and we determine the boundary conditions under the assumption that the fluid does not cross the boundary. Next we write a precise 
formulation of the problem to be studied. We then perform a change of dependent variables which reduces the original problem to one with velocity tangent to the boundary and then, by a change of independent variables, we transform the new equations into a system of PDE on a cylindrical domain.  Next we study a viscous regularization of our problem, showing
existence of weak solutions to the vorticity form of the Navier-Stokes system in a time-dependent domain with homogeneous Dirichlet boundary conditions. We then use the family of approximate solutions obtained and we 
obtain an priori estimate, which together with a compactness argument, enables us to choose a subsequence of the approximate solutions converging strongly in $L^2$ to a weak solution to the original problem.

Technically, the new feature of the present work is the treatment based the vorticity equation. The setup of the
problem is closely based on the original Navier-Stokes work, specially \cite{teramoto}. The existence result for
the vorticity equation on a noncylindrical domain with vorticity vanishing at the boundary is new, but its proof
is a standard proof based on Galerkin approximations. The a priori $L^p$ estimates independent of viscosity
are the main point of the argument. Finally, the passage to the limit is based on a standard compactness argument.

The main difference between our work and \cite{he} is the need for the first change of variables, which changes the
problem into a new one with velocity tangent to the moving boundary. This new change of variables introduces an extra 
convection term in the PDE of the form $\rho \nabla u$, where $\rho$ is a smooth, divergence-free vector field determined
by the boundary motion. The field $\rho$ is not tangent to the boundary. As a consequence, we cannot obtain an
energy estimate uniform in viscosity for the approximate problem, something that was needed in \cite{he}.

 The remainder of this work is divided into  four sections. In Section 2, we introduce basic notation and formulate
the material boundary condition. In Section 3 we perform the successive changes of variable which reduce the
problem to a system of PDE on a fixed domain. In Section 4 we construct an approximation of our problem based
on a viscous regularization. In Section 5 we prove a priori estimates and prove our main result.   

\section{Boundary conditions for moving domains}

In this section we discuss the boundary condition associated with a moving material boundary.
Let $Q_T = \bigcup_{0\leq t \leq T} \Omega_t \times \left\{t \right\}$ be a noncylindrical space-time domain, each $\Omega_t$ being a bounded domain in $\mathbb{R}^2$ with smooth boundary $\partial \Omega_t$.  In $Q_T$ we consider the initial-boundary value problem for the Euler equations:
 
\begin{equation} \label{seq1}
\left\{
\begin{array}{lcl}
\gd{\frac{\partial u}{\partial t}} + (u \cdot \nabla)u = -\nabla p + f && x \in \Omega_t, t>0\\\\
\mbox{div }u = 0 && x \in \Omega_t, t>0\\\\
u.\eta = g(x,t) && x \in \partial\Omega_t, t>0\\\\
u(x,0)=u_0(x) && x \in \Omega_0.
\end{array}
\right.
\end{equation}

\vspace{0.3cm}

Here $u=u(x,t)=(u_1,u_2)$ denote the unknown velocity and $p=p(x,t)$ denote the pressure of the ideal fluid at point $(x,t) \in \Omega_t \times \left\{t \right\}$, while $u_0(x)$ and \linebreak $f=f(x,t)$ denote respectively the given initial velocity vector field, and the external force vector field; $\eta$ is the unit outward normal vector of $\partial\Omega_t$ and $g=g(x,t)$ is given on the boundary $\bigcup_{0\leq t \leq T} \partial\Omega_t \times \left\{t \right\}$ from a prescribed movement of domain, assuming that the fluid does not cross the boundary.  We will describe how $g$ is
determined from $\Omega_t$ later in this section.

We impose the following conditions on the motion of the domain:

\begin{enumerate}
\item[(A.1)] $\overline{\Omega}_t$ is diffeomorphic to $\overline{\Omega}_0$ for each $t \in [0,T]$.
\item[(A.2)] The area of $\Omega_t$ is equal to that of $\Omega_0$ for each $t \in [0,T]$.
\end{enumerate}
 
As $\Omega_t$ has the smooth boundary $\partial\Omega_t$ and it moves smoothly with respect to t, we have the following lemma:

\begin{lemma} \label{lemma1}:Let $\gamma = \gamma(x,t)$ be the signed distance function with relation to $\Omega_t$.  Then there exists an open subset $U_t$ which contains $\partial\Omega_t$ such that $\gamma$ and $\partial_t \gamma$  are $C^\infty$ \linebreak functions in $U_t$.  Furthermore, $\left|\nabla \gamma(\overline{x},t) \right|=1$ for each $\overline{x} \in \partial\Omega_t$ and \linebreak $\partial\Omega_t=\left\{x:\gamma(x,t)=0\right\}.$
\end{lemma}
For more details and a proof, see \cite{evans}. Under the assumptions (A.1) and (A.2), we have the following result.

\begin{lemma} \label{lemma2}: There exists a cylindrical domain $\tilde{Q}_T = \tilde{\Omega} \times \mathbb{R}$ and a 
time-preserving diffeomorphism $\Phi:\overline{Q}_T \rt \overline{\tilde{Q}}_T$,
\[
(y,s) = \Phi(x,t) = (\Phi_1(x,t),\Phi_2(x,t),t)
\]
such that
\[
J^{-1}(t) = \mbox{det } \gd\left[\frac{\partial \Phi_i(x,t)}{\partial x_j}\right]_{i,j=1,2}\equiv 1 \hspace{0.5cm} (x,t) \in \displaystyle{Q}_T.
\]
\end{lemma}

For the proof see T. Miyakawa and Y. Teramoto \cite{teramoto}.
Beyond (A.1) and (A.2), we must include an assumption on the regularity 
of the diffeomorphism $\Phi$. More precisely, the domain $\Omega_t$ must satisfy the 
following assumption.

\vspace{3mm}

(A.3) The derivatives $\partial \Phi_i / \partial x_j$ and $\partial \Phi_i / \partial t$  $(1\leq i,j \leq2)$ are continuous and bounded functions on $\overline{Q}_T$.

\vspace{3mm}

Consider $Q_T$ a domain satisfying the assumptions (A.1), (A.2) and (A.3).  Let $u(x,t)$ be divergence-free vector field, defined in a neighborhood of $Q_T$, that makes $Q_T$ a material domain.  We have: 

\begin{proposition} \label{prop1}: There exists a unique $g(x,t)$ defined in $\bigcup_{0 \leq t \leq T} \partial \Omega_t \times \left\{t \right\}$ such that $u(x,t). \eta=g(x,t)$ in $\partial \Omega_t$, where $\eta:=\eta_t(x)$ is an unit outward normal vector in $x \in \partial\Omega_t$.  Furthermore, for each $t \in [0,T]$, we have \[\int_{\partial\Omega_t}g(x,t)ds=0.\]
\end{proposition}

{\it Proof:} We denote the inverse transformation of $\Phi(x,t)=(y,t)$ in Lemma \ref{lemma2} by $\Psi(y,t)=(\Psi_1(y,t),\Psi_2(y,t),t)$ where $(y_1,y_2) \in \tilde{\Omega}$ and $t \in [0,T]$. For $(\overline{x},t) \in \bigcup_{0 \leq t \leq T} \partial \Omega_t \times \left\{t \right\}$, set $(\overline{x},t)=(\Psi_1(\overline{y},t),\Psi_2(\overline{y},t),t)$ and by the Lemma \ref{lemma1} we have 

\[
\gamma(\Psi_1(\overline{y},t),\Psi_2(\overline{y},t),t)=0.
\]

Hence,

\begin{equation} \label{seqseq}
\frac{\partial \gamma(\overline{x},t)}{\partial t}+ \displaystyle\sum_{j=1}^{\infty}\frac{\partial \gamma(\overline{x},t)}{\partial x_j} \frac{\partial \Psi_j(\overline{y},t)}{\partial t}=0.
\end{equation}

If we denote $V_t(x)=(\displaystyle\frac{\partial \Psi_1(y,t)}{\partial t}, \displaystyle\frac{\partial \Psi_2(y,t)}{\partial t})$ where $(y,t)=\Phi(x,t)$, identity \eqref{seqseq} gives us

\[
V_t(\overline{x}).\eta_t(\overline{x})=g(\overline{x},t) \hspace{0.5cm}, \overline{x} \in \partial \Omega_t
\]
where 
\begin{equation}\label{seq2}
g(\overline{x},t)=\displaystyle\frac{\partial \gamma(\overline{x},t)}{\partial t}
\end{equation}
and 
\begin{equation}\label{seq3}
\eta_t(\overline{x})=-\nabla_x \gamma(\overline{x},t)
\end{equation}
by the Lemma \ref{lemma1}.

\vspace{0.3cm}

On the other hand, we can show that the vector field $V_t(x)$ is divergence-free on each $\Omega_t$ by using the Lemma \ref{lemma2}. Noticing that $\eta_t(\overline{x})$ is a unit normal vector at point $\overline{x} \in \partial\Omega_t$, we have by the divergence theorem:

\[
\int_{\partial \Omega_t} g(\overline{x},t) ds = \int_{\partial \Omega_t} V_t(\overline{x}).\eta_t(\overline{x}) ds =\int_{\Omega_t} \mbox{div }V_t(x)dx=0.
\]

Now, let $u(x,t)$ be divergence-free velocity vector field arbitrary, defined in a neighborhood of $Q_T$ that makes $Q_T$ a material domain.  The hypothesis that $Q_T$ is a material domain with respect to $u$, it allows to introduce the flow $X=X(\alpha,t)$ satisfying the ordinary differential equation in the introduction, and with the same argument used previously, we have for $g$ and $\eta_t$ as in (\ref{seq2}) and (\ref{seq3})
\[
u(\overline{x}).\eta_t(\overline{x})=g(\overline{x},t) \hspace{0.5cm}, \overline{x} \in \partial \Omega_t
\]
where $u(x,t)=\left(\displaystyle\frac{\partial X_1(y,t)}{\partial t},\displaystyle\frac{\partial X_2(y,t)}{\partial t}\right)$ for $(y,t)=X^{-1}(x,t).$  Which conclude this proof.

\section{Reduction to a fixed domain problem}

In this section we perform two changes of variables in order to reduce our problem to one suitable for
analytical treatment. The first one is a change of dependent variables designed to make the boundary 
condition homogeneous and the second one used the diffeomorphism $\Phi$ to change our problem to a
fixed-domain one. We begin with the homogenization of the boundary data.

In the formulation of the problem (\ref{seq1}), the unknown velocity field $u(.,t)$ has given normal component at each point in $\partial \Omega_t$ for each $t$.  We would like to transform the original problem to that of finding a velocity field without normal component.  To accomplish that, we first, consider $h$ the solution of the following problem.

\begin{equation} \label{seq4}
\left\{
\begin{array}{lcl}
\Delta h=0 &\mbox{em}&\Omega_t\\\\
\displaystyle\frac{\partial h}{\partial \eta}=g(.,t)&\mbox{em}&\partial \Omega_t.
\end{array}
\right.
\end{equation}

Since the Neumann boundary value $g(.,t)$ satisfies the compatibility condition $\int_{\partial \Omega_t} g(x,t)ds =0$ by Proposition \ref{prop1}, we can assume the existence of a $h \in C^{\infty}(\overline{\Omega}_t)$ satisfying (\ref{seq4}) for each $t$. For more details see \cite{strauss} and \cite{zac}.

Denote by $\rho \equiv \nabla h$ and observe that $\rho$ is a two-dimensional $C^{\infty}$ vector field on 
$\overline{Q}_T$ such that
\begin{equation}\label{seq5}
\begin{array}{lcl}
\mbox{div }\rho=0  &\mbox{em }& \Omega_t\\\\
\rho . \eta=g(x,t) &\mbox{em }&  \partial\Omega_t.
\end{array}
\end{equation}

We set
\[v(x,t)=u(x,t)-\rho(x,t)\]
in (\ref{seq1}), and we have the following equations:
\begin{equation} \label{seq6}
\left\{
\begin{array}{lcl}
\gd{\frac{\partial v}{\partial t}} + (v \cdot \nabla)v + (\rho \cdot \nabla)v+\nabla p =  f && x \in \Omega_t, t>0\\\\
\mbox{div }v = 0 && x \in \Omega_t, t>0\\\\
v.\eta = 0 && x \in \partial\Omega_t, t>0\\\\
v(x,0)=v_0(x) && x \in \Omega_0.
\end{array}
\right.
\end{equation}
where $p$ is now redefined as $p+ \partial_t h  + (1/2)(\rho \cdot \nabla)h + (1/2)(v \cdot \nabla)h$ and $v_0=u_0-\rho(.,0).$  Henceforth we will discuss the solvability of Equations (\ref{seq6}) for $\{v,p\}.$

Next we turn to the reduction to a cylindrical domain. Let
\begin{equation}\label{seq7}
\begin{array}{l}
\tilde{v}^i(y,s)=\displaystyle\frac{\partial y_i}{\partial x_1}v^1(\Phi^{-1}(y,s))+\displaystyle\frac{\partial y_i}{\partial x_2}v^2(\Phi^{-1}(y,s))\\\\
\tilde{\rho}^i(y,s)=\displaystyle\frac{\partial y_i}{\partial x_1}\rho^1(\Phi^{-1}(y,s))+\displaystyle\frac{\partial y_i}{\partial x_2}\rho^2(\Phi^{-1}(y,s))\\\\
\tilde{f}^i(y,s)=\displaystyle\frac{\partial y_i}{\partial x_1} f^1(\Phi^{-1}(y,s))+\displaystyle\frac{\partial y_i}{\partial x_2} f^2(\Phi^{-1}(y,s))\\\\
\tilde{p}(y,s)=p(\Phi^{-1}(y,s))
\end{array}
\end{equation}
for $i=1,2.$  Then (\ref{seq6}) is transformed into the following problem on $\tilde{Q}_T$ for \linebreak $\tilde{v}=(\tilde{v}_1,\tilde{v}_2)$ and $\tilde{p}$:

\begin{equation} \label{seq8}
\left\{
\begin{array}{lcl}
\gd{\frac{\partial \tilde{v}}{\partial s}} + M\tilde{v} + N_1\tilde{v}+ N_2\tilde{v}= \tilde{f} -\nabla_q \tilde{p} && y \in \tilde{\Omega}, s>0\\\\
\mbox{div }\tilde{v} = 0 && y \in \tilde{\Omega}, s>0\\\\
\tilde{v}.\tilde{\eta} = 0 && y \in \partial\tilde{\Omega}, s>0\\\\
\tilde{v}(y,0)=\tilde{v}_0(y) && y \in \tilde{\Omega}.
\end{array}
\right.
\end{equation}

\vspace{0.3cm}
Where $\tilde{v}_0=(\tilde{v}_0^1,\tilde{v}_0^2)$, $\tilde{\eta}$ denotes the unit exterior normal along $\partial\tilde{\Omega}$ and
$$
\begin{array}{l}
(M\tilde{v})=(\displaystyle\frac{\partial y_j}{\partial t})\nabla_j \tilde{v}_i+(\displaystyle\frac{\partial y_i}{\partial x_k})(\displaystyle\frac{\partial^2 x_k}{\partial s \partial y_j})\tilde{v}_j,\\\\
(N_1\tilde{v})^i=\tilde{\rho}_j\nabla_j \tilde{v}_i+\tilde{v}_j\nabla_j \tilde{\rho}_i, \hspace{0.8cm} (N_2\tilde{v})^i=\tilde{v}_j\nabla_j \tilde{v}_i,\\\\
q^{ij}=(\displaystyle\frac{\partial y_i}{\partial x_k})(\displaystyle\frac{\partial y_j}{\partial x_k}),\hspace{0.8cm} q_{ij}=(\displaystyle\frac{\partial x_k}{\partial y_i})(\displaystyle\frac{\partial x_k}{\partial y_j}),\\\\
(\nabla_q \tilde{p})^i= q^{ij} \displaystyle\frac{\partial \tilde{p}}{\partial y_j}, \hspace{0.8cm} \nabla_j \tilde{v}_i=\displaystyle\frac{\partial \tilde{v}_i}{\partial x_l}+\tilde{v}_k (\displaystyle\frac{\partial y_i}{\partial x_l})(\displaystyle\frac{\partial^2 x_l}{\partial y_j \partial y_k}).
\end{array}
$$

From now on, we use the summation convention, that is take sum over repeated indices.  Moreover, we let $\tilde{v}$ denote the vector field on $\tilde{Q}_T$ obtained by the transformation $\tilde{v}^i(y,s)=\partial y_i / \partial x_k . v^k(\Phi^{-1}(y,s))$ for each vector field $v$ on $Q_T$.  Conversely, $v$ is the vector field obtained by inverse transformation for $\tilde{v}$.  From the Lemma \ref{lemma2} it is easy to see that
\[
(q^{ij})^{-1}=(q_{ij}), \hspace{0.5cm} \mbox{det }(q_{ij})=J(t)^2.
\] 

We point out that the divergence operator is left invariant under the coordinate transformation.  Finally we note that $\partial \tilde{v} / \partial s +M\tilde{v}$, $N_1 \tilde{v}$ and $N_2 \tilde{v}$ correspond respectively to $\partial v / \partial t$, $(\rho \cdot \nabla) v+(v \cdot \nabla)\rho$ and $(v \cdot \nabla)v$ under the transformation $\Phi$; see \cite{inoue} for the details.

\begin{lemma}\label{lemma3}:

(1)The matrixes $\left[q_{ij}\right]$ and $\left[q^{ij}\right]$ are positive definite and bounded.

(2)The derivatives $\partial x_i/\partial y_j$ for i,j=1,2 are bounded functions on $\tilde{Q}_T.$
\end{lemma}

{\it Proof:} (1) If we denote $T=\left[\partial y_i/\partial x_k \right]$, then $\left[q^{ij}\right]=TT^t$ and we conclude that $\left<q^{ij}x,x\right>=\left\|T^tx\right\|^2\geq0$.  By the assumption $(A.3)$ the matrixes are bounded.
(2) It follows from the fact that $\left[\displaystyle\frac{\partial x_i}{\partial y_j}\right]=\left[\displaystyle\frac{\partial y_i}{\partial x_j}\right]^{-1}$, Lemma \ref{lemma2} and the assumption (A.3).

\vspace{0.3cm}
To show the existence of a weak solution for the system \eqref{seq6}, we will make use the solutions to the Navier-Stokes equations in a time dependent domain with a modified boundary condition. The next section will be dedicated to show existence of a weak solution to this approximate problem.

\section{The Navier-Stokes equations in a noncylindrical domain}

The purpose of this section is to construct a family of approximations which will be used to prove
existence of weak solutions for the ideal flow equations by means of a limit process. The approximate
problem will be the Navier-Stokes system, with a kind of slip boundary condition which is well-known to
behave well under vanishing viscosity. More precisely, we consider the system
  
\begin{equation} \label{seq9}
\left\{
\begin{array}{lcl}
\gd{\frac{\partial v}{\partial t}}-\nu \Delta v + (v \cdot \nabla)v + (\rho \cdot \nabla)v+\nabla p =  f && x \in \Omega_t, t>0\\\\
\mbox{div }v = 0 && x \in \Omega_t, t>0\\\\
v.\eta = 0, \omega=0 && x \in \partial\Omega_t, t>0\\\\
v(x,0)=v_0(x) && x \in \Omega_0.
\end{array}
\right.
\end{equation}	 
where $\omega=$curl $v$.

Our objective in this section is to prove the existence of a weak solution to problem \eqref{seq9}.
We begin by reducing the equations in (\ref{seq9}) to those in a cylindrical domain as in the subsection 
(3.2).  This yields:

\begin{equation} \label{seq10}
\left\{
\begin{array}{lcl}
\gd{\frac{\partial \tilde{v}}{\partial s}}-\nu L\tilde{v} + M\tilde{v} + N_1\tilde{v}+ N_2\tilde{v}= \tilde{f} -\nabla_q \tilde{p} && y \in \tilde{\Omega}, s>0\\\\
\mbox{div }\tilde{v} = 0 && y \in \tilde{\Omega}, s>0\\\\
\tilde{v}.\tilde{\eta} = 0,\gd{\left(\frac{\partial x_2}{\partial y_j}\frac{\partial y_i}{\partial x_1}-\frac{\partial x_1}{\partial y_j}\frac{\partial y_i}{\partial x_2}\right)\frac{\partial \tilde{v}^j}{\partial y_i}}=0 && y \in \partial\tilde{\Omega}, s>0\\\\
\tilde{v}(y,0)=\tilde{v}_0(y) && y \in \tilde{\Omega},
\end{array}
\right.
\end{equation}
where

\begin{equation}\label{seq*}
\omega=\left(\frac{\partial x_2}{\partial y_j}\frac{\partial y_i}{\partial x_1}-\frac{\partial x_1}{\partial y_j}\frac{\partial y_i}{\partial x_2}\right)\frac{\partial \tilde{v}^j}{\partial y_i}:=A:D_y\tilde{v},
\end{equation}

and

\vspace{0.3cm}

$(L\tilde{v})^i=q^{jk}\nabla_j \nabla_k \tilde{v}^i$, where

\vspace{0.3cm}

$\nabla_k \nabla_j \tilde{v}^i=\displaystyle\frac{\partial(\nabla_j \tilde{v}^i)}{\partial y_k}+\Gamma_{kl}^{i}\nabla_j \tilde{v}^l-\Gamma_{kj}^{l}\nabla_l \tilde{v}^i$ and

\vspace{0.3cm}

$\Gamma_{ij}^{k}=\left(\displaystyle\frac{\partial y_k}{\partial x_l}\right)\left(\displaystyle\frac{\partial^2 x_l}{\partial y_i \partial y_j}\right).$

\vspace{0.3cm}

The notation $L\tilde{v}$ correspond to $\Delta v$ under the transformation $\Phi$.

\vspace{0.3cm}

We introduce notation for some Hilbert spaces and inner products.  We denote by $\tilde{H}$ the space 
of square-integrable vector fields $(L^2(\tilde{\Omega}))^2$ and by $\tilde{V}$ the space $(H^1_0(\tilde{\Omega}))^2$. Similarly, we define 
$H_t$ and $V_t$ as the same spaces based on the domain $\Omega_t$. For each $t \in \R$, 
$H^1(\tilde{\Omega})$ is a Hilbert space with respect to the inner product:
\begin{equation}\label{seq11}
\left<\tilde{u},\tilde{v}\right>_t= \int_{\tilde{\Omega}}q_{ij}(y,t)\tilde{u}_i(y)\tilde{v}_j(y)dy,
\end{equation}
for $u,v \in H^1(\tilde{\Omega})$.

In $H_t$ we consider the usual inner product
\begin{equation}\label{seq12}
\left(u,v\right)_t= \int_{\Omega_t}u(x) \cdot v(x)dx.
\end{equation}

For $\tilde{u},\tilde{v} \in \tilde{V}$ the inner product in $\tilde{V}$ is denoted by

\begin{equation}\label{seq13}
\left<\left< \tilde{u},\tilde{v} \right> \right>_t \equiv \left<\nabla_q \tilde{u},\nabla_q \tilde{v}\right>_t= \int_{\tilde{\Omega}}q_{ij}(y,t)q^{kl}(y,t) \nabla_k\tilde{u}_i(y)\nabla_l\tilde{v}_j(y)dy.
\end{equation}

We consider the change of variables $x = \Phi^{-1}(y,t)$ and 
\begin{equation} \label{gab}
\tilde{u}^i(y,s)=\displaystyle\frac{\partial y_i}{\partial x_1}u^1(\Phi^{-1}(y,s))+\displaystyle\frac{\partial y_i}{\partial x_2}u^2(\Phi^{-1}(y,s)),
\end{equation}
with the same relation between $\tilde{v}$ and $v$. Note that, under this
change of variables, for any fixed $t$, (\ref{seq11}) is transformed into (\ref{seq12}) and 
(\ref{seq13}) is transformed into

\begin{equation}\label{seq14}
\left(\left( u,v\right) \right)_t= \int_{\Omega_t}\partial_{x_j}u(x) \cdot \partial_{x_j}v(x)dx.
\end{equation}

The norms corresponding to the inner products (\ref{seq12}) are denoted by $||\cdot||_t$.

Next we introduce  $\omega = \mbox{ curl }v$ and $w = \mbox{ curl }u$, and we introduce 

\begin{equation}\label{seq15}
\begin{array}{l}
\omega=\left(\displaystyle\frac{\partial x_2}{\partial y_j}\displaystyle\frac{\partial y_i}{\partial x_1}-\displaystyle\frac{\partial x_1}{\partial y_j}\displaystyle\frac{\partial y_i}{\partial x_2}\right)\displaystyle\frac{\partial \tilde{v}^j}{\partial y_i}:=A:D_y\tilde{v}= A:D_yK_{\tilde{\Omega}}\left[\tilde{\omega}\right]\\\\
w=\left(\displaystyle\frac{\partial x_2}{\partial y_l}\displaystyle\frac{\partial y_k}{\partial x_1}-\displaystyle\frac{\partial x_1}{\partial y_l}\displaystyle\frac{\partial y_k}{\partial x_2}\right)\displaystyle\frac{\partial \tilde{u}^l}{\partial y_k}:=B:D_y\tilde{u}= B:D_yK_{\tilde{\Omega}}\left[\tilde{u}\right]
\end{array}
\end{equation}
where $\tilde{v}=K_{\tilde{\Omega}}\left[\tilde{\omega}\right]$ and $\tilde{u}=K_{\tilde{\Omega}}\left[\tilde{w}\right]$ are given by the {\it Biot-Savart} law since div $\tilde{v}=0$ and div $\tilde{u}=0.$

Observe that

\begin{eqnarray*}
\frac{\partial \omega}{\partial x_l}& = &\left(\frac{\partial^2y_i}{\partial x_l \partial x_1}\frac{\partial x_2}{\partial y_j}-\frac{\partial^2y_i}{\partial x_l \partial x_2}\frac{\partial x_1}{\partial y_j}\right)\frac{\partial \tilde{v}^j}{\partial y_i}+A\left(\frac{\partial^2 \tilde{v}^j}{\partial y_i \partial y_k}\right)\frac{\partial y_k}{\partial x_l}:=\\
&:=&E_1:D_y K_{\tilde{\Omega}}\left[\tilde{\omega}\right]+A D_{y}^{2}K_{\tilde{\Omega}}\left[\tilde{\omega}\right]\frac{\partial y_{k_{1}}}{\partial x_l},
\end{eqnarray*}
in the same way
\begin{eqnarray*}
\frac{\partial w}{\partial x_l}:=E_2:D_y K_{\tilde{\Omega}}\left[\tilde{w}\right]+B D_{y}^{2}K_{\tilde{\Omega}}\left[\tilde{w}\right]\frac{\partial y_{k_{2}}}{\partial x_l}.
\end{eqnarray*}

We define an alternative inner product in $L^2(\tilde{\Omega})$ as follows. 

\begin{equation} \label{seq16}
\left[\tilde{\omega},\tilde{w}\right]_{t}=\int_{\tilde{\Omega}}\left(A:D_yK_{\tilde{\Omega}}\left[\tilde{\omega}\right]\right)\left(B:D_yK_{\tilde{\Omega}}\left[\tilde{w}\right]\right)dy,
\end{equation}
and $\left(\omega,w\right)_t$ is the usual inner product in $L^{2}(\Omega_t)$.  

Since we obtain $(\ref{seq16})$ from $\left(\omega,w\right)_t$ by implementing the usual coordinate transformation, 
we can conclude that $(\ref{seq16})$ is, in fact, an inner product.
Finally, $H_{\sigma}^{2}(\Omega_t)$ is the subspace of divergence-free vector fields in $H^2(\Omega_t)$. 

Next we introduce a convenient notion of weak solution for the initial-boundary
value problem (\ref{seq9}).

\begin{definition}\label{def1}
A velocity field $v \in L^2(0,T;H_{\sigma}^{2}(\Omega_t))\cap L^{\infty}(0,T;V_t)$, for any $T>0$ being fixed, is called a {\it weak solution} of system (\ref{seq9}) with initial data $v_0(x)$ and forcing $f(x,t)$, if for any $\theta(.,t)\in C^1_c([0,T);V_t)$, the following identities are satisfied:

\begin{enumerate}
\item[i)]
\[
-\int_0^T (v(t),\theta^{'}(t))_t dt+\nu \int_0^T ((v(t),\theta(t)))_t dt+\int_0^T ((v(t) \cdot\nabla) v(t),\theta(t))_t dt+
\]
\[
+\int_0^T ((\rho \cdot\nabla) v(t),\theta(t))_t dt=(v_0,\theta(\cdot,0))_0+\int_0^T (f(t),\theta(t))_t dt,
\]

\item[ii)] the velocity is incompressible in the weak sense, that is 
\[
\int_0^T\int_{\Omega_t}\nabla \theta(t)v(t)dx dt=0,
\]

\item[iii)]$v.\eta=0$ in $\partial \Omega_t \times (0,T),$

\vspace{0.5cm}

\item[iv)] $\omega=0$ in $\partial \Omega_t \times (0,T).$

\end{enumerate}
\end{definition}

{\bf Remark:} The space of test functions $C^1_c([0,T),V_t)$ may be defined by the condition that 
$\tilde{u} \in C^1_c([0,T),V_t)$ if and only if the $u$ associated with $\tilde{u}$ through \eqref{gab}
belongs to $C^1_c([0,T),\tilde{V})$.

Now, we state and prove existence of a weak solution to (\ref{seq9}).

\begin{theorem}\label{teo1}
Fix an arbitrary $T>0$.  Then for each $v_0 \in V_0$, $\mbox{curl }v_0 \in L^2(\Omega_0)$ and each $f \in L^2(0,T;V_t)$ there exists a weak solution of (\ref{seq9}) satisfying Definition \ref{def1}.
\end{theorem}

We will obtain a weak solution to the problem (\ref{seq9}) through its vorticity formulation.  
At the level of strong solutions, if $v = v(x,t)$ is a solution of \eqref{seq9}, then the associated
vorticity $\omega = \mbox{ curl }v$, we have the following equations

\begin{equation} \label{seq17}
\left\{
\begin{array}{lcl}
\gd{\frac{\partial \omega}{\partial t}}-\nu \Delta \omega + (v \cdot \nabla)\omega + (\rho \cdot \nabla)\omega = \mbox{curl } f && \Omega_t, t>0\\\\
\omega=0 && \partial\Omega_t, t>0\\\\
\omega(x,0)=\omega_0(x) && \Omega_0.
\end{array}
\right.
\end{equation}	 

\vspace{0.3cm}
Moreover, taking $\tilde{\omega}=$ curl $\tilde{v}$ and applying the curl in the first equation of the system (\ref{seq17}), we get an expression for that first equation defined in $\tilde{\Omega}$ as follows.
\[
\frac{\partial \tilde{\omega}}{\partial s}-\nu(\mbox{curl}(L \tilde{v}))+\mbox{curl}(M\tilde{v})+\mbox{curl} (N_1\tilde{v}+N_2\tilde{v})=\mbox{curl}(\tilde{f}).
\]

\vspace{0.3cm}
We will prove Theorem \ref{teo1} by constructing a family of approximate solutions using Galerkin. 
The construction will proceed as follows. Let $\{\tilde{\zeta}_j\}$ be a sequence of linearly independent vectors in $C_{0}^{\infty}(\tilde{\Omega})$ total in $H_{0}^{1}(\tilde{\Omega})$, and $\{\tilde{\alpha}_j(y,t)\}$ be its Schmidt orthogonalization with respect to the inner product (\ref{seq16}).  Note that $\tilde{\alpha}_j(t)=\tilde{\alpha}_j(\cdot,t)$ thus obtained is smooth in $(\cdot,t)$, because it is a finite linear combination of $\{\tilde{\zeta}_j\}$  with coefficients in $C^{\infty}([0,T];\mathbb{R}).$

Taking $\tilde{\theta}_j = K_{\tilde{\Omega}}[\tilde{\alpha}_j]$, there exists a stream function $\tilde{\psi}_j$ such that $\tilde{\theta}_j=\nabla^{\perp}\tilde{\psi}_j$, and  $\tilde{\psi}_j$ satisfies:

$$
\left\{
\begin{array}{ccl}
\Delta \tilde{\psi}_j=\tilde{\alpha}_j &\mbox{in}& \tilde{\Omega}\\\\
\tilde{\psi}_j=0 &\mbox{in}& \partial\tilde{\Omega}.
\end{array}
\right.
$$
 
We define approximate solutions $\tilde{\omega}_m(t)$, $m \in \mathbb{N}$, by the following equations
\begin{equation}\label{seq18}
\tilde{\omega}_m(t)=\sum_{j=1}^{m}h_{jm}(t)\tilde{\alpha}_j(t),
\end{equation}

\begin{equation}\label{seq19}
\tilde{\omega}_m(0)=\sum_{j=1}^{m}h_{jm}^0(t)\tilde{\alpha}_j(0), \hspace{0.5cm} h_{jm}^0=\left[ \tilde{\omega}_0,\tilde{\alpha}_j(0) \right]_0,
\end{equation}
where $\{h_{jm}(t)\}$ with $1 \leq j\leq m$, is defined by

\begin{equation}\label{seq20}
\left[\tilde{\omega}_m^{'},\tilde{\alpha}_j\right]_t=\left[\mbox{ curl }(\nu L \tilde{v}_m-M\tilde{v}_m-N_1\tilde{v}_m-N_2\tilde{v}_m+\tilde{f}),\tilde{\alpha}_j\right]_t.
\end{equation}

\vspace{0.3cm}

It is easy to see that for each $m \in \mathbb{N}$, there exists $t^*=t^*(m)$ such that $\tilde{\omega}_m(t)$ is determined uniquely by (\ref{seq20}) if $1 \leq j\leq m$ and $t \in[0,t^*]$, as defined in (\ref{seq18}) satisfying (\ref{seq19}).  The next lemma guarantees that $\tilde{\omega}_m(t)$ is defined on the whole interval $[0,T].$

\begin{lemma}\label{lemma4}:
$\{\omega_m(t)\}$ remains bounded in $L^{\infty}(0,T;L^2(\Omega_t))\cap L^{2}(0,T;H^1(\Omega_t))$.
\end{lemma}
{\it Proof:}  We rewrite (\ref{seq20}) in $\Omega_t$ to obtain:
\begin{equation}\label{seq**}
(\omega_m^{'},\alpha_j)_t-(\nu \Delta\omega_m,\alpha_j)_t=-(v_m \cdot \nabla\omega_m+\rho \cdot \nabla\omega_m,\alpha_j)_t+(\mbox{curl }f,\alpha_j).
\end{equation}

Multiply it by $h_{jm}(t)$ and take the sum in $j$ to get
\begin{equation}\label{seq21}
(\omega_m^{'},\omega_m)_t-(\nu \Delta\omega_m,\omega_m)_t=-(v_m \cdot \nabla\omega_m+\rho \cdot \nabla\omega_m,\omega_m)_t+(\mbox{curl }f,\omega_m)_t.
\end{equation}

\vspace{0.3cm}

Integrating by parts and using that $\omega_m$ has compact support in $\Omega_t$ we have
\[
\frac{1}{2}\frac{d}{dt}\left|\left|\omega_m\right|\right|_t^2+\nu \int_{\Omega_t}  \left|\nabla \omega_m\right|^2 dx=-\int_{\Omega_t}\frac{\partial v_m^i}{\partial x_i}\left|\omega_m\right|^2 +\frac{\partial \rho^i}{\partial x_i}\left|\omega_m\right|^2 dx+(\mbox{curl }f,\omega_m)_t.
\]

\vspace{0.3cm}

By the Holder and Young inequalities, remembering that div $v_m =$ div $\rho=0$, and integrating in $t$, we obtain
\begin{equation}\label{seq22}
\left|\left|\omega_m\right|\right|_t^2+2 \nu \int_0^t\left|\left|\nabla\omega_m\right|\right|_{\sigma}^2d\sigma \leq \left|\left|\omega_0\right|\right|_0^2+\int_0^t\left|\left|\mbox{curl }f\right|\right|_{\sigma}^2d\sigma+ \int_0^t\left|\left|\omega_m\right|\right|_{\sigma}^2d\sigma. 
\end{equation}

\vspace{0.3cm}

By Gronwall's lemma;
\[
\left|\left|\omega_m\right|\right|_t^2 \leq (\left|\left|\omega_0\right|\right|_0^2+\int_0^T \left|\left|\mbox{curl }f\right|\right|_{t}^2dt)(1+Te^T).
\]
From that inequality and (\ref{seq22}) we conclude this proof.

In order to prove Theorem \ref{teo1} we require compactness of the approximating sequence $\{\omega_m\}$.
This compactness follows from the {\it a priori} estimates in Lemma \ref{lemma4} by means of an argument
which is a straightforward adaptation of a similar result, see Lemma 2.5 in \cite{teramoto}. 
We choose not to repeat this argument here, and we will just state the corresponding fact as a lemma, 
ommiting the proof.   

\begin{lemma}\label{lemma5}:
$\{\omega_m(t)\}$ is precompact in $L^{2}(0,T;L^2(\Omega_t))$.
\end{lemma}
 
We are now ready to prove Theorem \ref{teo1}.

{\bf Proof of Theorem \ref{teo1}:}

By Lemma \ref{lemma4} and Lemma \ref{lemma5} we may assume, passing to subsequence as necessary, that there exists 
$\omega \in L^{\infty}(0,T;L^2(\Omega_t))\cap L^{2}(0,T;H^1(\Omega_t))$ such that $\omega_m$ converges to $\omega$ in $L^{2}(0,T;H^1(\Omega_t))$ weakly, $\omega_m$ converges to $\omega$ in $L^{\infty}(0,T;L^2(\Omega_t))$ weak-star and $\omega_m$ converges to $\omega$ in $L^{2}(0,T;L^2(\Omega_t))$ strongly.

Now we integrate (\ref{seq**}) in $t$ and integrate by parts to get:
\begin{eqnarray*}
-\int_0^T\int_{\Omega_t}\omega_m\frac{\partial \alpha_j}{\partial t}dx dt +\nu \int_0^T\int_{\Omega_t} \nabla\alpha_j \nabla\omega_m dx dt+\int_0^T\int_{\Omega_t}\alpha_j(v_m\cdot\nabla)\omega_m dx dt+\\
+\int_0^T\int_{\Omega_t}\alpha_j(\rho\cdot\nabla)\omega_m dx dt=\int_{\Omega_0}\omega_m(0)\alpha_j(x,0)dx+\int_0^T\int_{\Omega_t}\alpha_j(\mbox{curl }f)dx dt.
\end{eqnarray*}

Since $\omega_m(0) \rt \omega_0$ in $L^2(\Omega_0)$ and $v_m$ is bounded in $L^{\infty}(0,T;H^1(\Omega_t))$, by letting $m \rt \infty$ we obtain

\begin{eqnarray*}
-\int_0^T\int_{\Omega_t}\omega\frac{\partial \alpha_j}{\partial t}dx dt +\nu \int_0^T\int_{\Omega_t} \nabla\alpha_j \nabla\omega dx dt+\int_0^T\int_{\Omega_t}\alpha_j(v\cdot\nabla)\omega dx dt+\\
+\int_0^T\int_{\Omega_t}\alpha_j(\rho\cdot\nabla)\omega dx dt=\int_{\Omega_0}\omega_0\alpha_j(x,0)dx+\int_0^T\int_{\Omega_t}\alpha_j(\mbox{curl }f)dx dt.
\end{eqnarray*}

By linearity this equality holds for $\alpha=\sum_{j=1}^{l}\alpha_j$. Recall that $\{\alpha_j\}$ is total in $H_{0}^{1}(\Omega_t)$.  Therefore, for any $\alpha \in H_{0}^{1}(\Omega_t)$ the previous equality holds.  
Since $\alpha$ is arbitrary, we consider $\alpha = \psi(t)$ where $\psi(t)$ is the stream function associated
with the flow $\theta$. We have that $\theta=\nabla^{\perp}\psi$, and therefore,
$$
\left\{
\begin{array}{ccc}
\Delta \psi=\mbox{curl }\theta &\mbox{in }& \Omega_t\\
\psi=0 &\mbox{in }& \partial\Omega_t
\end{array}
\right.
$$
and $\theta= K_{\Omega_t}[\alpha].$  Since $\omega=$ curl $v$ we obtain

\begin{eqnarray*}
-\int_0^T\int_{\Omega_t}(\mbox{curl }v)\frac{\partial \psi(t)}{\partial t}dx dt +\nu \int_0^T\int_{\Omega_t}\nabla(\mbox{curl }v) \nabla(\psi(t))dx dt+\\
+\int_0^T\int_{\Omega_t} (\psi(t))(v\cdot\nabla)(\mbox{curl }v) dx dt+\int_0^T\int_{\Omega_t}(\psi(t))(\rho\cdot\nabla)(\mbox{curl }v) dx dt=\\
=\int_{\Omega_0}(\mbox{curl }v_0(x))(\psi(0))(x,0)dx+\int_0^T\int_{\Omega_t}(\psi(t))(\mbox{curl }f)dx dt.
\end{eqnarray*}

Integrating by parts and using the fact that $\theta=\nabla^{\perp}\psi$ we get,

\begin{eqnarray*}
-\int_0^T\int_{\Omega_t}v\frac{\partial \theta}{\partial t}dx dt +\nu \int_0^T\int_{\Omega_t}\nabla v \nabla \theta dx dt+\int_0^T\int_{\Omega_t} \theta(v\cdot\nabla) v dx dt+\\
+\int_0^T\int_{\Omega_t}\theta(\rho\cdot\nabla)v dx dt=\int_{\Omega_0}v_0(x)\theta(x,0)dx+\int_0^T\int_{\Omega_t}\theta f dx dt.
\end{eqnarray*}
which satisfies the first condition on the Definition \ref{def1}.  The other conditions are proved by integration 
by parts and taking traces in a straightforward manner.  This concludes the proof.

\section{Euler equations in a noncylindrical domain}

Now we are ready to show a result of existence of weak solution to the Euler equations (\ref{seq6}) defined in a time dependent domain. 

\begin{definition}\label{def2} Let $T>0$. A velocity field $v \in L^{\infty}(0,T;L^2(\Omega_t))$ is called a {\it weak solution} of the Euler equations (\ref{seq6}) with initial data $v_0(x)$ and external force field $f(x,t)$, if for any $\tilde{\varrho}=h(t)\tilde{\theta}$ such that $\tilde{\theta}\in \tilde{V}$ and $h \in C^1([0,T],\mathbb{R})$, $h(T)=0$, the following identities are satisfied:
\begin{enumerate}
\item[i)]
\[
-\int_0^T \left<\tilde{v}(t),\tilde{\varrho}^{'}(t)\right>_t dt-\int_0^T \left<\tilde{v}(t) ,M\tilde{\varrho}(t)\right>_t dt+
\]
\[
+\int_0^T \left<N_1\tilde{v}(t)+N_2\tilde{v}(t),\tilde{\varrho}(t)\right>_t dt=\left<\tilde{v}_0,\tilde{\varrho}(0)\right>_0+\int_0^T \left<\tilde{f}(t),\tilde{\varrho}(t)\right>_t dt,
\]

\item[ii)] the velocity is incompressible in the weak sense, that is 
\[
\int_0^T\int_{\tilde{\Omega}}\nabla \tilde{\varrho}(t)\tilde{v}(t)dy dt=0,
\]

\item[iii)]$\tilde{v}.\tilde{\eta}=0$ in $\partial \tilde{\Omega} \times (0,T).$
\end{enumerate}
\end{definition}

We will prove existence of a weak solution to the Euler equations in the sense of Definition \ref{def2}. 
We will assume that the external force field $f \in L^2(0,T;W^{1,r}(\Omega_t))$ for $1<r\leq \infty$ is potential, 
that is, $f=\nabla\overline{f}$ where $\overline{f} \in L^2(0,T;W^{2,r}(\Omega_t))$ for $1<r\leq \infty$.

\begin{theorem}\label{teo2}
Fix $1<r\leq \infty$. For $v_0 \in V_0$, with $\mbox{curl } v_0 \in L^r(\Omega_0)$ and $f=\nabla\overline{f}$ where $\overline{f} \in L^2(0,T;W^{2,r}(\Omega_t))$, there exists a weak solution $v \in L^{\infty}(0,T;L^2(\Omega_t))$ of (\ref{seq6}) satisfying
$$
\begin{array}{lcl}
\nabla v, \omega \in L^{\infty}(0,T;L^r(\Omega_t))&if&1<r<\infty\\\\
\nabla v \in L^{\infty}(0,T;L^p(\Omega_t)),\omega \in L^{\infty}(0,T;L^{\infty}(\Omega_t))&if&1<p<\infty \hspace{0.2cm}for \hspace{0.2cm} r=\infty,
\end{array}
$$
where $\omega=\mbox{curl }v.$
\end{theorem}

We use Theorem \ref{teo1}to construct an approximate solution sequence. For $\omega_0 \in L^r(\Omega_0)$,
we consider a sequence $\omega_0^{\nu}$ of smooth functions defined in $\Omega_0$ such that
$\omega_0^{\nu} \to \omega_0$ strongly in $L^r(\Omega_0)$ when $\nu \to 0$. Such a sequence can be
obtained, for example, by solving the heat equation with homogeneous Dirichlet conditions in $\Omega_0$ with
initial data $\omega_0$ for time $\nu$. For fixed $T,\nu>0$, we use Theorem \ref{teo1} to obtain $\omega_{\nu}\in L^{\infty}(0,T;L^{2}(\Omega_t)\cap L^{2}(0,T;H^1(\Omega_t))$, weak solution to the system:

\begin{equation} \label{seq23}
\left\{
\begin{array}{lcl}
\gd{\frac{\partial \omega_{\nu}}{\partial t}}-\nu \Delta \omega_{\nu} + (v_{\nu} \cdot \nabla)\omega_{\nu} + (\rho \cdot \nabla)\omega_{\nu} = 0 && \Omega_t, t>0\\\\
\omega_{\nu}=0 && \partial\Omega_t, t>0\\\\
\omega_{\nu}(x,0)=\omega^{\nu}_0(x) && \Omega_0.
\end{array}
\right.
\end{equation}	

In particular, for any $\alpha \in H_{0}^{1}(\Omega_t)$, we have:

\begin{equation}\label{seq24}
(\omega_{\nu}^{'},\alpha)_t=(\nu \Delta \omega_{\nu},\alpha)_t-(v_{\nu}\cdot\nabla \omega_{\nu},\alpha)_t-(\rho\cdot\nabla \omega_{\nu},\alpha)_t.
\end{equation}

\vspace{0.3cm}
We look for an {\it a priori} estimate for the vorticity $\omega_{\nu}$, uniform in $\nu$. Let $\phi_{\epsilon}:\mathbb{R}\rt \mathbb{R}$ be defined as follows

$$
\phi_{\epsilon}(x)=\left\{
\begin{array}{ccl}
|x|^r &\mbox{if }&|x|\geq \epsilon \\
\frac{r}{2}\epsilon^{r-2}x^2+\epsilon^{r}(1-r/2)&\mbox{if }&\left|x\right|\leq \epsilon.
\end{array}
\right.
$$

\vspace{0.3cm}

Observe that $\phi_{\epsilon}$ is a $C^1$ and convex function, and $\phi_{\epsilon}(0)=0.$  
We have the following result

\begin{proposition}\label{teo3}
For any $t \in [0,T]$, $\{\omega_{\nu}(t)\}$ is bounded in $L^{\infty}(0,T;L^r(\Omega_t))$ \linebreak for $1<r\leq\infty.$
\end{proposition}

{\it Proof:} The key issue is that $\phi_{\epsilon}^{'}(\omega_{\nu}) \in H_0^1(\Omega_t)$, 
and therefore we can use it as a test function in the definition of weak solution.
We can rewrite (\ref{seq24}) as follows:
$$
\begin{array}{c}
\displaystyle\int_{\Omega_t}\displaystyle\frac{d}{dt}\phi_{\epsilon}(\omega_{\nu})dx=\nu \displaystyle\int_{\Omega_t}( \Delta \phi_{\epsilon}(\omega_{\nu})
-\phi_{\epsilon}^{''}(\omega_{\nu})\left|\nabla \omega_{\nu}\right|^2)dx-\displaystyle\int_{\Omega_t} \mbox{div }(v_{\nu}\phi_{\epsilon}(\omega_{\nu}))dx-\\\\
-\displaystyle\int_{\Omega_t}\mbox{div }(\rho\phi_{\epsilon}(\omega_{\nu}))dx=\nu \displaystyle\int_{\partial\Omega_t}\nabla\phi_{\epsilon}(\omega_{\nu}).\eta ds-\nu \displaystyle\int_{\Omega_t}\phi_{\epsilon}^{''}(\omega_{\nu})\left|\nabla \omega_{\nu}\right|^2 dx-\\\\
-\displaystyle\int_{\partial\Omega_t}\phi_{\epsilon}(\omega_{\nu})v_{\nu} \cdot \eta ds-\displaystyle\int_{\partial\Omega_t}\phi_{\epsilon}(\omega_{\nu})\rho \cdot \eta ds=\nu\displaystyle\int_{\partial\Omega_t}\phi_{\epsilon}^{'}(0)\nabla\omega_{\nu}\cdot\eta ds-\\\\
-\nu \displaystyle\int_{\Omega_t}\phi_{\epsilon}^{''}(\omega_{\nu})\left|\nabla \omega_{\nu}\right|^2dx-\phi_{\epsilon}(0)\displaystyle\int_{\partial\Omega_t}(v_{\nu}.\eta+g(x,t))ds\leq0.
\end{array}
$$

Integrating this inequality in $t$ we have:

\begin{equation}\label{seq25}
\int_{\Omega_t}\phi_{\epsilon}(\omega_{\nu}(t))dx \leq \int_{\Omega_0}\phi_{\epsilon}(\omega_{\nu}(0))dx.
\end{equation}

As $\{\phi_{\epsilon}(x)\}$ is a increasing sequence and $\displaystyle\lim_{\epsilon \to 0} \phi_{\epsilon}(x) =\left|\omega_{\nu}\right|^r$, we have by the Fatou's Lemma that
\[
\int_{\Omega_t}\left|\omega_{\nu}(t)\right|^r \leq \displaystyle\lim_{\epsilon \to 0} \mbox{inf}\int_{\Omega_t}\phi_{\epsilon}(\omega_{\nu}(t))dx.
\]

Therefore, if $1<r<2$ we obtain
\begin{equation}\label{seq26}
\left|\left|\omega_{\nu}\right|\right|_{L^r(\Omega_t)}\leq \left|\left| \omega_0^{\nu}\right|\right|_{L^r(\Omega_0)}
\leq C, 
\end{equation}
because converging sequences are bounded. If $2\leq r<\infty$ it is enough to define $\phi(x)=\left|x\right|^r$, 
to conclude that $\phi^{'}(\omega_{\nu}) \in H_0^1(\Omega_t)$ and to proceed in the same way as in the case $1<r<2$. 
Therefore we get (\ref{seq26}) for $1<r\leq\infty$, which completes the proof.

{\it Proof of Theorem \ref{teo2}:}

By the Calder\'on-Zygmund theorem and Proposition \ref{teo3}, we have the following estimate
\begin{equation}\label{seq27}
\left|\left|\nabla v_{\nu} \right|\right|_{L^r(\Omega_t)}\leq C \left|\left|\omega_{\nu} \right|\right|_{L^r(\Omega_t)}\leq C \left|\left|\mbox{curl }v_0\right|\right|_{L^r(\Omega_0)},
\end{equation}
for $1<r<\infty.$

Moreover, if $r>1$,
\[
C \leq C_1 \frac{r^2}{r-1}
\]
where $C_1$ is a constant independent of $r$.  This fact follows by tracking the constant in 
the Marcinkiewicz interpolation inequality, see \cite{gilbarg}.

As
\[
\tilde{\omega}_{\nu}=\mbox{curl }\tilde{v}_{\nu}=\frac{\partial x_i}{\partial y_1}\frac{\partial y_2}{\partial x_j}\frac{\partial v_{\nu}^j}{\partial x_i}-\frac{\partial x_i}{\partial y_2}\frac{\partial y_1}{\partial x_j}\frac{\partial v_{\nu}^j}{\partial x_i},
\]
we deduce, with the help of the assumption $(A.3)$, Lemma \ref{lemma3} and estimate (\ref{seq27}), that
\begin{equation}\label{seq28}
\left|\left|\tilde{\omega}_{\nu} \right|\right|_{L^r(\tilde{\Omega})}\leq \mbox{C }.
\end{equation}

By using the Calder\'on-Zygmund theorem once more we get

\begin{equation}\label{seq29}
\left|\left|\nabla\tilde{v}_{\nu} \right|\right|_{L^r(\tilde{\Omega})}\leq \mbox{C } \left|\left|\tilde{\omega}_{\nu} \right|\right|_{L^r(\tilde{\Omega})}\leq \mbox{C },
\end{equation}
for $1<r<\infty.$

From the inequalities (\ref{seq26}), (\ref{seq27}) and (\ref{seq29}) it follows that there exists a subsequence of $\{\omega_{\nu}\}$ (without relabeling) such that
\begin{equation}\label{seq30}
\left\{
\begin{array}{cccc}
\omega_{\nu} \stackrel{*}{\fraco} \omega &\mbox{in }&L^{\infty}(0,T;L^r(\Omega_t))&1<r\leq \infty\\
\nabla v _{\nu} \stackrel{*}{\fraco}\nabla v &\mbox{in }&L^{\infty}(0,T;L^r(\Omega_t))&1<r<\infty\\
\nabla \tilde{v} _{\nu} \stackrel{*}{\fraco}\nabla \tilde{v} &\mbox{in }&L^{\infty}(0,T;L^r(\tilde{\Omega}))&1<r<\infty.\\
\end{array}
\right.
\end{equation}

Moreover, observation 2.6 on \cite{teramoto} shows that $\{\tilde{v}_{\nu}\}$ is precompact on $L^{\infty}(0,T;L^2(\tilde{\Omega}))$, and this guarantees the existence of $\tilde{v} \in L^{\infty}(0,T;L^2(\tilde{\Omega}))$ such that
\begin{equation}\label{seq311}
\tilde{v}_{\nu} \rt \tilde{v} \hspace{0.3cm} \mbox{in } \hspace{0.3cm} L^{\infty}(0,T;L^2(\tilde{\Omega})).
\end{equation}

As $\omega^{\nu}$ satisfies Definition \ref{def2} we have that for any $\alpha \in H_0^1(\Omega_t)$ the following
identity holds
\begin{equation}\label{seq31}
\begin{array}{c}
-\displaystyle\int_0^T \displaystyle\int_{\Omega_t}\omega_{\nu}\displaystyle\frac{\partial \alpha}{\partial t}dx dt+\nu \int_0^T \int_{\Omega_t}\nabla\omega_{\nu} \nabla \alpha dx dt+\displaystyle\int_0^T \displaystyle\int_{\Omega_t}\alpha(v_{\nu}\cdot\nabla)\omega_{\nu}dx dt+\\\\
+\displaystyle\int_0^T \displaystyle\int_{\Omega_t}\alpha(\rho\cdot\nabla)\omega_{\nu}dx dt=\displaystyle\int_{\Omega_0}\omega_{\nu}(x,0)\alpha(x,0)dx+\displaystyle\int_0^T \displaystyle\int_{\Omega_t}\alpha(\mbox{curl }f)dx dt.
\end{array}
\end{equation}

We consider $\alpha=h(t)\psi$ where $h \in C^1([0,T];\mathbb{R})$ with $h(T)=0$ and $\psi$ is understood as a
stream function with respect to a test velocity field $\theta$. We have that $\theta=\nabla^{\perp}\psi$, and
$$
\left\{
\begin{array}{ccc}
\Delta \psi=\mbox{curl }\theta &\mbox{in }& \Omega_t\\\\
\psi=0 &\mbox{in }& \partial\Omega_t.
\end{array}
\right.
$$

Since $\omega_{\nu}=$ rot $v_{\nu}$, we integrate by parts and we use the fact that $\theta=\nabla^{\perp}\psi$ to obtain

\[ -\displaystyle\int_0^T \displaystyle\int_{\Omega_t}v_{\nu}\displaystyle\frac{\partial(h(t)\theta)}{\partial t}dx dt
+ \nu \int_0^T \int_{\Omega_t}\nabla v_{\nu} \nabla (h(t)\theta) dx dt\]
\[ +\displaystyle\int_0^T \displaystyle\int_{\Omega_t}(h(t)\theta)(v_{\nu}\cdot\nabla)v_{\nu}dx dt
+\displaystyle\int_0^T \displaystyle\int_{\Omega_t}(h(t)\theta)(\rho\cdot\nabla)v_{\nu}dx dt\]
\[=\displaystyle\int_{\Omega_0}v_{\nu}(x,0)(h(0)\theta(x,0))dx+\displaystyle\int_0^T \displaystyle\int_{\Omega_t}(h(t)\theta)f dx dt.\]
 
Set $\varrho=h(t)\theta$ and rewrite this identity in $\tilde{Q}_T$ to obtain

$$
\begin{array}{c}
-\displaystyle\int_0^T \left<\tilde{v}_{\nu},\tilde{\varrho}^{'}(t)\right>_t dt-\displaystyle\int_0^T \left<\tilde{v}_{\nu},M\tilde{\varrho}(t)\right>_t dt+\nu \displaystyle\int_0^T \left<\nabla_q\tilde{v}_{\nu},\nabla_q\tilde{\varrho}(t)\right>_t dt+\\\\
+\displaystyle\int_0^T \left<N_1\tilde{v}_{\nu}+N_2\tilde{v}_{\nu},\tilde{\varrho}(t)\right>_t dt=\left<\tilde{v}_{\nu}(x,0),\tilde{\varrho}(0)\right>_0+\displaystyle\int_0^T \left<\tilde{f},\tilde{\varrho}(t)\right>_t dt.
\end{array}
$$

We let $\nu \rt 0$, to prove thatthe first condition of the Definition \ref{def2} is satisfied.

As the norm is weakly lower-semicontinuous, (\ref{seq26}) and (\ref{seq27}) give us the estimates
$\left|\left|\omega(t)\right|\right|_{L^r(\Omega_t)}\leq \left|\left|\omega_{\nu}(t)\right|\right|_{L^r(\Omega_t)}<\infty$ (for $t>0$ and $1<r\leq\infty$) and $\left|\left|\nabla v(t)\right|\right|_{L^r(\Omega_t)}\leq \left|\left|\nabla v_{\nu}(t)\right|\right|_{L^r(\Omega_t)}<\infty$ (for $t>0$ and $1<r<\infty$).

If $r=\infty$, then $\omega_0 \in L^p(\Omega_0)$ when $1<p<\infty$.  Therefore (\ref{seq27}) gives us the estimate
$\left|\left|\nabla v(t)\right|\right|_{L^r(\Omega_t)}<\infty$ when $t>0$ and $1<p<\infty$ if $r=\infty$,
which concludes the proof of Theorem \ref{teo2}.

We would like to conclude this paper by mentioning a couple of natural questions that arise naturally 
from our analysis. First, is it possible to generalize this result to include 
initial vorticities in $L^1$, or nonnegative bounded measures, extending Delort's Theorem to 
noncylindrical domains? We refer the reader to \cite{delort} and \cite{vecchiwu} for the relevant 
existence results. We have proved our existence theorem relying only on vorticity estimates
because we do not have an useful energy estimate in this context. After
multiplying by $u$ and integrating by parts, the term $\rho \cdot \nabla v$ in \eqref{seq9} gives rise
to a boundary term of the form:
\[ \int_{\partial \Omega_t} \frac{|v|^2}{2} \rho \cdot \eta dS =  \int_{\partial \Omega_t} \frac{|v|^2}{2} g(x,t)dS,\]
and we cannot control such a boundary term in a manner that is independent of viscosity. This absence of an 
energy estimate is the key technical difference between our work and He and Hsiao's. Now, Delort's Theorem, and
its adaptation to $L^1$ vorticities by Vecchi and Wu, require in an essential manner {\it a priori} estimates both
for vorticity and kinetic energy, and the observation above makes this extension a difficult problem.
It would probably be technically challenging but doable to extend Delort's Theorem to He and Hsiao's context.    
The second natural question is whether one can remove the condition that $\Omega_t$ be simply connected.
One special case would be the flow on the exterior of a moving body. This would be interesting from
the physical point of view. Since our approach is based on vorticity estimates, extending our result to
domains with holes depends on understanding and controlling the harmonic part of the flow for noncylindrical
domains.  

\small{ {\it Acknowledgments:} The authors would like to thank J.-L. Boldrini for pointing out 
a serious error in an early version of this work and H. J. Nussenzveig Lopes for useful conversations
and comments. This work is based on F. Z. Fernandes Ph.D. dissertation, done under the supervision of 
M. C. Lopes Filho, in connection with the graduate program in mathematics at the State University of 
Campinas (UNICAMP). The graduate work of F. Z. Fernandes was supported by FAPESP grant \# 01/06486-5; 
M. C. Lopes Filho's research is supported  in part by CNPq grants \# 302.102/2004-3 and \# 472504/2004-5. }

\end{document}